\documentclass[12pt,leqno]{amsart}

\usepackage{enumerate} 
\usepackage{booktabs} 
\usepackage[utf8]{inputenc} 
\usepackage[T1]{fontenc} 
\usepackage{amssymb,amsthm}
\usepackage{mathrsfs}
\usepackage{textcomp}
\usepackage{bm}
\usepackage{graphicx}
\usepackage{xcolor}

\newtheorem{thm}{Theorem}[section]

\theoremstyle{definition}

\theoremstyle{remark}

\numberwithin{equation}{section}

\usepackage[font=footnotesize,tableposition=top,figureposition=bottom,skip=0pt]{caption}
\usepackage{hyperref}
\hypersetup{%
	linkcolor=blue,
	anchorcolor=black,
	citecolor=red,
	filecolor=magenta,
	menucolor=red,
	urlcolor=magenta,
	colorlinks=true,
	bookmarksnumbered=true,
	bookmarksopen=true,
	bookmarksopenlevel=0,
	breaklinks=true
}

\begin{document}

\title[Oscillating sequences and Multiple ergodic limits]{Fully oscillating sequences and weighted multiple ergodic limit\footnote{Note in CRAS Paris}}

\author{Ai-hua Fan}
\address{LAMFA, UMR 7352 CNRS, University of Picardie, 33 rue Saint Leu,80039 Amiens, France}
\email{ai-hua.fan@u-picardie.fr}
\thanks{}
\maketitle






\begin{abstract} We  prove that
fully oscillating sequences  are orthogonal to multiple ergodic realizations of affine  maps of zero entropy on compact abelian groups.
It is more than what Sarnak's conjecture requires for these dynamical systems. 
\end{abstract}

\section{Introduction and results}
A
     sequence of complex numbers $ c = (c_n)_{n\ge 0} $ is said to be \textit{oscillating of order $ d $} ($ d\ge 1 $) 
	 if for any real polynomial $ P\in \mathbb{R}_d[x] $  of degree less than or equal to $ d $ we have
	 \begin{equation}\label{eq:OSCd}
	 \lim\limits_{N\to \infty}\frac{1}{N}\sum_{n=0}^{N-1}c_{n}e^{2\pi i{P}(n)} = 0.
	 \end{equation}
	 It is said to be {\em fully oscillating}  if it is oscillating of all orders.
	 This notion of oscillation of higher order was introduced in \cite{F}. The oscillation of order $1$ was earlier considered in \cite{FJ} in order to formulate some results related to Sarnak's conjecture.
	 See \cite{Sarnak}, \cite{Sarnak2} for Sarnak's conjecture. See \cite{ALR}, \cite{AKL}, \cite{DG},  \cite{FKLM}, \cite{HLSY}, \cite{KL}, \cite{HWZ},  \cite{Wang} for some  related recent works.
	 The M\"{o}bius sequence $(\mu(n))$ is a typical example of  fully oscillating sequences (\cite{Davenport}, \cite{Hua}). Recall that $\mu(1)=1$,  $\mu(n)=(-1)^k$ if $n$ is square free and has $k$ distinct prime factors, and $\mu(n)=0$ for other integers $n$.  Random subnormal sequence is almost surely fully oscillating (\cite{F}) and the sequence $(\beta^n \mod 1)$
	 is fully oscillating for almost all $\beta >1$ (\cite{AJ}).
	 The oscillating sequences of orders $d$ are characterized by their orthogonality to different classes of dynamical systems (\cite{Shi}).  
	 
	 Sarnak's conjecture states that for any topological dynamical system $(X, T)$ of zero entropy,  for any continuous function $f \in C(X)$ and any point $x \in X$, we have
	 $$
	        \lim_{N\to \infty} \frac{1}{N}\sum_{n=1}^N \mu(n) f(T^n x) =0.
	 $$
	 This conjecture remains open in its generality. The above equality is referred to as the orthogonality of M\"{o}bius sequence to the realization $f(T^n x)$ of the system $(X, T)$, or as 
	 the disjointness of $(\mu(n))$ to the system $(X, T)$.
	 For $\ell (\ge 1)$ functions $f_1, \cdots, f_\ell \in C(X)$, the sequence $f_1(T^n x) f_2(T^{2n} x) \cdots f_\ell(T^{\ell n} x)$ could be referred to as a multiple ergodic realization.
	 \medskip

Following Liu and Sarnak \cite{LS}, we can prove the following orthogonality of fully oscillating sequences to the multiple ergodic realizations of affine linear maps on a compact abelian group which are of zero entropy. 	
\medskip 
	\begin{thm}\label{thm:1}
		Let $ \ell\ge 1 $ be an integer. Let $G$ be a compact abelian group. Assume that\\	
		 \indent \textnormal{(i)} \ \  $T: G \to G$ is an affine linear map of zero entropy. \\
	 \indent \textnormal{(ii)}\  $ (c_n) $ is a fully oscillating sequence.\\
	 \indent \textnormal{(iii)} $ q_1,\cdots,q_{\ell}\in \mathbb{Z}[x] $ are $ \ell $ polynomials such that $ 
	 q_j(\mathbb{N})\subset 
	\mathbb{N} $ for all $j$.  
	 \newline
	Then for any continuous function $ F \in C(G^\ell) $ and any point $ x\in G $, we have
	\begin{equation}\label{eq:1-1}
	\lim\limits_{N\to\infty}\dfrac{1}{N}\sum_{n=0}^{N-1}c_n F(T^{q_1(n)}x, \cdots, T^{q_\ell(n)}x) = 0.
	\end{equation}
	\end{thm}
	
	Recall that an affine  map on a compact abelian group $G$
	  is defined by
	  $$
	         T x = A x +b
	  $$ 
	  where $A: G \to G$ is an automorphism of $G$ and $b \in G$.
	
	Theorem \ref{thm:1} generalizes the following theorem due to Liu and Sarnak which holds for M\"{o}bius sequence to fully oscillating sequences. 
	\medskip
	
	\begin{thm} [Liu and Sarnak \cite{LS}] \label{ThmLS}
	 The M\"{o}bius sequence $(\mu(n))$ is  linearly disjoint from any affine linear map on a  compact abelian group which is of zero entropy.
	\end{thm} 
	\medskip
	
	The result of Theorem \ref{thm:1} was proved in \cite{F}, based on \cite{HP0}, \cite{HP1}, \cite{HP2}, for the class of topological systems of quasi-discrete spectrum in the sense of Hahn-Parry \cite{HP0}, including minimal affine linear maps on a connected compact abelian group.
	The proof of Theorem \ref{thm:1} in this note will be based on ideas of Liu and Sarnak and on the fact that arithmetic subsequences of oscillating sequences are oscillating.  One of the ideas of Liu and Sarnak
	is stated as follows. It is drawn from the proof of their first theorem in \cite{LS}.
	\medskip
	
\begin{thm} [Liu-Sarnak \cite{LS}]\label{LemmaLS}
Let $Tx = Ax +b$ be an affine linear map of zero entropy on $X:=\mathbb{T}^d\times F$ where $d\ge 1$ and $F$ is a finite abelian group.
Consider the automorphism $W$ on the product group $X\times X$ defined by $W(x_1, x_2):= (Ax_1 + x_2, x_2)$. Then
\\
\indent {\rm (i)} \ \ $(T^n x, b) = W^n(x, b)$ for all $x \in X$.\\
\indent {\rm (ii)} \ there exist integers $\nu \ge 1$ and $\kappa \ge 0$ such that $W^\nu = I  + N $ where $N$ is nilpotent in the sense $N^{\kappa+1}=0$. 
\end{thm}
\medskip

The following fact which has its own interest will also be useful in the proof of Theorem \ref{thm:1}.
\medskip
	
\begin{thm} \label{PO}
Let $a\ge 2$ be an integer. 
A sequence  $(w_n)$ is oscillating of order $d$ if and only if the arithmetic subsequences $(w_{an+b})$ are all oscillating of order $d$ for any integer $b\ge 0$.
\end{thm}
\medskip

Before proving Theorem \ref{thm:1} we prove Theorem \ref{PO}.

\section{Proof of Theorem \ref{PO}}

If $(w_{an +b})$ are oscillating of order $d$ for $0\le b\le a-1$, it is obvious that 
$(w_n)$ is oscillating of order $d$. 

Now assume that $(w_n)$ is oscillating of order $d$. 
First observe that from the definition, it is clear that any shifted sequence
$(w_{n +b})$ ($b\ge 0$) is oscillating of order $d$. So, it suffices to prove that $(w_{an})$ is oscillating of order $d$. 
Since $a$ can be decomposed into product of primes, we have only to prove that  $(w_{p n})$ is oscillating of order $d$
for any prime  $p\ge 2$.

Let $P\in \mathbb{R}_d[t]$. For any $N\ge p$, denote 
\begin{eqnarray*}
      S_N  &=& \sum_{0\le n <N} w_n e^{2\pi i P(n)}, \\
      S_{N, j}  &=& \sum_{m: 0\le pm +j <N} w_{pm+j} e^{2\pi i P(pm +j)} \ \ \ (0\le j <p).
\end{eqnarray*}
We have the trivial decomposition 
\begin{equation}\label{RX}
   S_N = S_{N,0} + S_{N,1} +\cdots + S_{N, p-1}. 
\end{equation}
For any integer $0\le u\le p-1$, denote 
\begin{eqnarray*}
      S_N^u  &=& \sum_{0\le n <N} w_n e^{2\pi i (P(n) + \frac{u n}{p})}
\end{eqnarray*}
where $x \mapsto P(x) + u x/p$ is a real polynomial.
Write $\omega=e^{2\pi i /p}$.  We have the following decomposition
\begin{equation}\label{RX2}
   S_N^u  = S_{N,0} + \omega^u S_{N,1} + \omega^{2u} S_{N, 2}\cdots + \omega^{(p-1)u}S_{N, p-1},
\end{equation}
which is similar to (\ref{RX}) and which contains (\ref{RX}) as a particular case corresponding to $u=0$.
Taking  sum over $u$, we get
$$
   \sum_{u=0}^{p-1}S_{N}^u = p S_{N,0} + S_{N,1}\sum_{u=0}^{p-1}\omega^{u} + S_{N,2}\sum_{u=0}^{p-1}\omega^{2 u} +\cdots +
   S_{N,p-1}\sum_{u=0}^{p-1}\omega^{(p-1)u}.
$$
Since $p$ is prime, any $j$ with  $1\le j \le p-1$ is invertible in the ring $\mathbb{Z}/p\mathbb{Z}$ so that $\sum_{u=0}^{p-1}\omega^{j u} =0$ for all $1\le j \le p-1$. Thus
$$
        S_{N, 0} = \frac{1}{p} \sum_{u=1}^{p-1} S_{N}^u.
$$
Since the sequence $\{w_n\}$ is oscillating of order $d$, we have $S_N^u = o(N)$ for all $u$ as $N\to \infty$, so that 
$$
    \lim_{N\to \infty} \frac{1}{N} \sum_{m=0}^{[N/p]} w_{pm} e^{2\pi i P(pm)} =   \lim_{N\to \infty} \frac{S_{N, 0}}{N}=0. 
$$
This implies that $m \mapsto w_{pm}$ is oscillating of order $d$.

\section{Proof of Theorem \ref{thm:1}}

Let $\widehat{G}$ be the dual group of $G$. We have $\widehat{G^\ell} = \widehat{G}^\ell$.  Any continuous function $F \in C(G^\ell)$ can be uniformly approximated by
trigonometric polynomials on $G^\ell$, which are finite linear combinations of functions of the form $\phi_1(x_1) \cdots \phi_\ell(x_\ell)$ where $\phi_1, \cdots, \phi_\ell \in \widehat{G}$. So,  it suffices to prove that
	\begin{equation}\label{eq:1-2}
	\lim\limits_{N\to\infty}\dfrac{1}{N}\sum_{n=0}^{N-1}c_n \phi_1(T^{q_1(n)}x) \cdots \phi_\ell( T^{q_\ell(n)}x) = 0
	\end{equation}
	holds for all $\phi_1, \cdots, \phi_\ell \in \widehat{G}$ and all $x\in G$ (see \cite{F} for details). 
	
	Now, we mimick Liu and Sarnak \cite{LS}.  First remark that the problem can be reduced to torus. In fact, let $\Phi =\{\phi_1, \cdots, \phi_\ell\} \subset \widehat{G}$.
	Recall that the action of $A$ on $\widehat{G}$ is defined by $(A \gamma) (x) = \gamma(Ax)$ for $\gamma \in \widehat{G}$ and $x\in G$. 
         Let $\langle \Phi\rangle$ be the smallest $A$-invariant closed subgroup of $\widehat{G}$ which contains $\Phi$ and
         $$\langle \Phi\rangle^\perp:=\{x \in G: \gamma(x)=1\ \  \forall \gamma \in \widehat{G}\}$$ 
         be the annihilator of $\langle \Phi\rangle$, a closed subgroup of $G$. Let $$G_\Phi: =G/\langle \Phi \rangle^\perp$$ be the quotient group.
         The $A$-invariance of $\langle \Phi \rangle$ implies that 
         $$
             \forall x\in G, \ \forall y \in \langle \Phi\rangle^\perp, \ \  T(x+y) = Tx \ \mod \langle \Phi\rangle^\perp.
         $$
        Thus $T$ induces an affine map,  which will be denoted by $T_\Phi$,  on the quotient group $G_\Phi$. Being a factor of $(G, T)$, the system
        $(G_\Phi, T_\Phi)$ has zero entropy. By Aoki's Theorem (\cite{Aoki}, a statement in the proof on p.13),  $\langle \Phi \rangle$ is finitely generated.
        Remark that $\widehat{G}_\Phi = \langle \Phi \rangle$. Then $G_\Phi$, as dual group of $\langle \Phi \rangle$, is isomorphic to $\mathbb{T}^d \times F$ for some $d\ge 1$ and some finite abelian group $F$. 
        On the other hand, for any $\phi\in \Phi$ and  any $x\in G$, we have
        $$
                   \phi(T^n x) = \widetilde{\phi} (T_\Phi^n \widetilde{x})
        $$
        where $ \widetilde{x} := x + \langle \Phi \rangle^\perp$ is the projection of $x$ onto $G_\Phi$ and   $\widetilde{\phi}$ is the character in $\widetilde{G}_\Phi$ induced by $\phi$.
        Thus the proof of (\ref{eq:1-2}) is reduced to the dynamics $(\mathbb{T}^d \times F, T_\Phi)$.
        
        By Theorem \ref{LemmaLS}, it suffices to treat the automorphism $W$ appearing in  Theorem \ref{LemmaLS}. Recall that $W^\nu = 1 + N$
        with $N^{\kappa+1} =0$. For any integer $n\ge 1$, write $n = m\nu + r$ with $m\ge 0$ and $0\le r \le \nu-1$.  The following expression was obtained in \cite{LS}
        \begin{equation}\label{E_LS}
              W^{m \nu + r} x= \sum_{j=0}^{\min (m, \kappa)} \binom{m}{j}
                                   N^j y_{r, j} =   \sum_{j=0}^{\kappa}  \binom{m}{j}   N^j y_{r, j}
        \end{equation}
        when $m \ge \kappa$, where $y_{r, j} = W^r N^j x$. Notice that these $y_{r, j}$ with $0\le r \le \nu-1$ and $0\le j \le \kappa$ are independent of $m$.  
        For any polynomial $q$ (a typical polynomial of $q_1, \cdots, q_\ell$), we have $q(m\nu +r) = q'(m)\nu +r'$
        where $r' = q(r) \mod \nu$  ($0\le r'\le \nu-1$) is independent of $m$  too, and $q' \in \mathbb{Z} [z]$ is a polynomial having the same  degree  as $q$. It follows from (\ref{E_LS}) that
         \begin{equation}\label{E_LS2}
              W^{q(m \nu + r)} x =W^{q'(m) \nu + r'} x=    \sum_{j=0}^{\kappa} \binom{q'(m) }{ j} 
                                    N^j y_{r', j}.
        \end{equation}
        This holds for $m$ sufficiently large so that $q'(m)\ge \kappa$. 
        Apply (\ref{E_LS2}) to each $q:=q_s$ ($1\le s \le \ell$). Then
        $$
              \phi_s(W^{q_s(m\nu +r)}) =  \prod_{j=0}^{\kappa} \phi_s(N^j y_{r_s', j})^{\binom{q_s'(m)}{j}} 
        $$
        Let $t_{s, r, j}$ be  the argument of the complex number $\phi_s(N^j y_{r_s', j})$. Then we get
        $$
              \prod_{s=1}^\ell \phi_s(W^{q_s(m\nu +r)})  = e^{2\pi i P(m)}
        $$
        where $P\in \mathbb{R}[x]$ is the real polynomial
        $$
              P(x) = \sum_{s=1}^\ell \sum_{j=0}^\kappa t_{s,r,j} \binom{q_s'(m)}{j}.
        $$
        By Theorem \ref{PO}, for any $0\le r <\nu$ we get 
        $$
         \lim_{M\to \infty}\frac{1}{M} \sum_{m=1}^M w_{m\nu +r} \prod_{s=1}^\ell \phi_s(W^{q_s(m\nu +r)} x) =
          \lim_{M\to \infty}\frac{1}{M} \sum_{m=1}^M w_{m\nu +r} e^{2\pi i P(m)}=0.
        $$
        Finally 
        $$
            \lim_{N\to \infty}\frac{1}{N} \sum_{n=0}^{N-1} w_{n} \prod_{s=1}^\ell \phi_s(W^{q_s(n)} x) =
            0.
        $$
        
        \medskip
        {\em Addendum.}  The oscillation of order $d$ is strongly related to the control of the $(d+1)$-th Gowers uniformity norm (see \cite{Tao}, \cite{Konieczny}). We can use Gowers uniformity norms
        to study oscillating properties. 


\begin{thebibliography}{99}

\bibitem{ALR} E. H. El Abdalaoui, M. Lema\'nczyk and T. de la Rue, 
{\em Automorphisms with quasi-discrete spectrum, multiplicative functions and average orthogonality along short intervals},
 International Mathematics Research Notices (2016), to appear.
 
 \bibitem{AKL}
 E. H. El Abdalaoui, S. Kasjan, and M. Lema\'nczyk, {\em 0-1 sequences of the Thue-Morse type and Sarnak's conjecture}, Proc. Amer. Math. Soc. 144 (2016), no. 1, 161-176.

\bibitem{AJ} S. Akiyama and Y. P. Jiang, {\em Higher order oscillation and uniform distribution}, preprint.


\bibitem{Aoki}
N. Aoki, {\em Topological entropy of distal affine transformations on compact abelian groups}, J. Math. Soc. Japan 23 (1971), 11-17.

\bibitem{Davenport}
H. Davenport, {\em On some infinite series involving arithmetical functions (II)},
Quart. J. Math. Oxford, 8 (1937), 313-320.

\bibitem{DG}
 T. Downarowicz and E. Glasner, {\em Isomorphic extensions and applications},
Topological Methods in Nonlinear Analysis (2015), to appear,



\bibitem{F} A. H. Fan, {\em Oscillating sequences of higher orders and topological systems of
quasi-discrete spectrum}, preprint.

\bibitem{FJ}
A. H. Fan and Y. P. Jiang, {\em Oscillating sequences, minimal mean attractability
and minimal mean-Lyapunov-stability}, Erg. Th. Dynam. Syst., to appear.

\bibitem{FKLM} S. Ferenzi, J. Kulaga-Przymus, M. Lema\'nczyk, and C. Mauduit, {\em Substitutions and M\"{o}bius disjointness}, preprint (2015).

\bibitem{HP0} F. Hahn and W. Parry,\,\emph{Minimal dynamical systems with quasi-discrete spectrum}. 
J. London Math. Soc., \textbf{40}\,(1965), 309-323.

\bibitem{HP1} H. Hoare and W. Parry,\,\emph{Affine transformations with quasi-discrete spectrum (I)}. 
J. London Math. Soc., \textbf{41}\,(1966), 88-96.

\bibitem{HP2} H. Hoare and W. Parry,\,\emph{Affine transformations with quasi-discrete spectrum (II)}. 
J. London Math. Soc., \textbf{41}\,(1966), 529-530.



\bibitem{Hua}
L. G. Hua, {\em Additive Theory of Prime Numbers} (Translations of Mathematical
Monographs : Vol 13). Amer Mathematical Society. 1966.

\bibitem{HLSY} W. Huang, Z. Lian, S. Shao and X. Ye, {\em Sequences from zero entropy noncommutative toral automorphisms and Sarnak Conjecture}, preprint.

\bibitem{HWZ}
W. Huang, Z. R. Wang and G. H. Zhang, {\em M\"{o}bius disjointness for topological models of ergodic systems with discrete spectrum},
preprint.

\bibitem{Konieczny} J. Konieczny, {\em Gowers norms for the Thue-Morse and Rudin-Schapiro sequences},  preprint

\bibitem{KL} J. Kulaga-Przymus and M. Lema\'nczyk, {\em The M\"{o}bius function and continuous extensions of rotations}, Monatsh. Math. 178 (2015), no. 4, 553-582.

\bibitem{LS}
 J. Y. Liu and P. Sarnak, {\em The M\"{o}bius function and distal flows}, Duke Math. J.
 164 (2015), no. 7, 1353-1399.
\bibitem{Sarnak}
P. Sarnak, {\em Three lectures on the M\"{o}bius function}, randomness and dynamics,
IAS Lecture Notes, 2009;\\
http://publications.ias.edu/sites/default/files/MobiusFunctionsLectures(2).pdf.
\bibitem{Sarnak2}
P. Sarnak, {\em M\"{o}bius randomness and dynamics}, Not. S. Afr. Math. Soc. 43
(2012), 89-97.

\bibitem{Shi} R. X. Shi, {\em Equivalent definitions of oscillating sequences of higher orders},
preprint.

\bibitem{Tao}T. Tao, {\em Higher order Fourier analysis}, volume 142 of Graduate Studies in Mathematics. American Mathematical Society, Providence, RI, 2012.

\bibitem{Wang} Z. R. Wang, {\em M\"{o}bius disjointness for analytic skew products}, preprint.

\end{thebibliography}
\medskip
{\em Acknowledgement.} This work is partially supported by NSFC 11471132. The author, supported by Knuth and Alice Wallenberg Foundation, visited Lund University in the autumn  2016 
where  the work was done.

\end{document}